\documentclass{article}

\usepackage{amsfonts, amsmath, amssymb, amscd, amsthm}
\usepackage{graphicx}

\newtheorem{thm}{Theorem}[section]
\newtheorem{cor}[thm]{Corollary}
\newtheorem{lem}[thm]{Lemma}

\theoremstyle{definition}
\newtheorem{defn}[thm]{Definition}
\theoremstyle{remark}

\newcommand{\G }{\Gamma (G, X\cup \mathcal H)}
\newcommand{\Ga }{\Gamma (G, \mathcal A)}
\newcommand{\dxh }{dist_{X\cup\mathcal H}}

\newcommand{\Hl }{\{ H_\lambda , \lambda \in \Lambda \} }
\newcommand{\e }{\varepsilon }

\newfont{\eufm}{eufm10}

\begin{document}

\title{Elementary subgroups of relatively hyperbolic groups and bounded generation.}
\author{D. V. Osin \thanks{This work has been supported by the RFFR Grant $\sharp $
02-01-00892.}}
\date{}%

\maketitle

\begin{abstract}
Let $G$ be a group hyperbolic relative to a collection of
subgroups $\{ H_\lambda ,\lambda \in \Lambda \} $. We say that a
subgroup $Q\le G$ is  hyperbolically embedded into $G$, if $G$ is
hyperbolic relative to $\{ H_\lambda ,\lambda \in \Lambda \} \cup
\{ Q\} $. In this paper we obtain a characterization of
hyperbolically embedded subgroups. In particular, we show that if
an element $g\in G$ has infinite order and is not conjugate to an
element of $H_\lambda $, $\lambda \in \Lambda $, then the (unique)
maximal elementary subgroup contained $g$ is hyperbolically
embedded into $G$. This allows to prove that if $G$ is boundedly
generated, then $G$ is elementary or $H_\lambda =G$ for some
$\lambda \in \Lambda $.
\end{abstract}

\section{Introduction}

Originally the notion of a relatively hyperbolic group was
proposed by Gromov in order to generalize various examples of
algebraic and geometric nature. Gromov's idea has been elaborated
by Bowditch \cite{Bow} in terms of the dynamics of group actions
on hyperbolic spaces, and by Farb \cite{Farb} in terms of the
geometry of Cayley graphs. Another definition of relative
hyperbolicity of a group $G$ with respect to a collection of
subgroups $\{ H_\lambda ,\lambda \in \Lambda \} $ is suggested in
\cite{Osin03}. In contrast to the definitions of Bowditch and
Farb, this approach does nor require the group $G$ and the
subgroups $H_\lambda $ to be finitely generated, as well as the
collection $\{ H_\lambda ,\lambda \in \Lambda \} $ to be finite.
This is important for some applications (see \cite{Osin04}) and,
in particular, allows to include the small cancellation theory
over free products developed in \cite[Ch. 5]{LS} within the
general frameworks of the theory of relatively hyperbolic groups.
On the other hand, in case the group $G$ is finitely generated our
definition is equivalent to the definitions of Bowditch and Farb
\cite{Osin03}.

More precisely, let $G$ be a group, $\Hl $ a collection of
subgroups of $G$, $X$ a subset of $G$. We say that $X$ is a {\it
relative generating set of $G$ with respect to $\Hl $} if $G$ is
generated by the set $\left(\bigcup\limits_{\lambda \in \Lambda}
H_\lambda \right) \cup X$. (We always assume that $X$ is
symmetrized, i.e. $X^{-1}=X$.) In this situation the group $G$ can
be regarded as the quotient group of the free product
\begin{equation}
F=\left( \ast _{\lambda\in \Lambda } H_\lambda  \right) \ast F(X),
\label{F}
\end{equation}
where $F(X)$ is the free group with the basis $X$. Let $N$ denote
the kernel of the natural homomorphism $F\to G$. If $N$ is a
normal closure of a finite subset $\mathcal R\subseteq N$ in the
group $F$, we say that $G$ has {\it relative presentation}
\begin{equation}\label{G}
\langle X,\; H_\lambda, \lambda\in \Lambda \; |\; R=1,\; R\in
\mathcal R \rangle .
\end{equation}
If $\sharp\, X<\infty $ and $\sharp\, \mathcal R<\infty $, the
relative presentation (\ref{G}) is called {\it finite} and the
group $G$ is called {\it finitely presented relative to the
collection of subgroups $\Hl $.}

Let
\begin{equation}\label{H}
\mathcal H=\bigsqcup\limits_{\lambda\in \Lambda} (H_\lambda
\setminus \{ 1\} )
\end{equation}
(we regard $H_\lambda $ as subgroups of $F$ here). Given a word
$W$ in the alphabet $X\cup \mathcal H$ such that $W$ represents
$1$ in $G$, there exists an expression
\begin{equation}
W=_F\prod\limits_{i=1}^k f_i^{-1}R_if_i \label{prod}
\end{equation}
with the equality in the group $F$, where $R_i\in \mathcal R$ and
$f_i\in F $ for any $i$. The smallest possible number $k$ in a
representation of type (\ref{prod}) is denoted by $Area^{rel}(W)$.

\begin{defn}\label{IP}
We say that a function $f:\mathbb N\to \mathbb N$ is a {\it
relative isoperimetric function} of (\ref{G}) if for any $n\in
\mathbb N$ and any word $W$ over  $X\cup \mathcal H$ of length $\|
W\| \le n$ representing the identity in the group $G$, we have
$Area^{rel} (W)\le f(n).$ The smallest relative isoperimetric
function of (\ref{G}) is called the {\it relative Dehn function}
of $G$ with respect to $\{ H_\lambda , \lambda \in \Lambda \} $
and is denoted by $\delta^{rel}_{G,\, \{ H_\lambda , \lambda \in
\Lambda \} }$ (or simply by $\delta ^{rel}$ when the group $G$ and
the collection of subgroups are fixed).
\end{defn}

We note that $\delta ^{rel} (n)$ is not always well--defined,
i.e., it can be infinite for certain values of the argument, since
the number of words of bounded relative length can be infinite.
Indeed consider the group
$$G=\langle a,b \; |\; [a,b]=1 \rangle\cong \mathbb Z\times
\mathbb Z $$ and the cyclic subgroup $H$ generated by $a$. Clearly
$X=\{ b^{\pm 1}\} $ is a relative generating set of $G$ with
respect to $H$. It is easy to see that the word $W_n=[a^n,b]$ has
length $4$ as a word over $X\cup H$ for every $n$, but $Area^{rel}
(W_n)$ growths linearly as $n\to\infty $. Thus we have $\delta
^{rel}(4)=\infty $ in this case.

However if $\delta ^{rel}$ is well--defined, it is independent of
the choice of the finite relative presentation up to the following
equivalence relation \cite[Theorem 2.32]{Osin03}. Two functions
$f,g:\mathbb N\to \mathbb N$ are called {\it equivalent} if there
are positive constants $A,B,C$ such that $f(n)\le Ag(Bn)+Cn$ and
$g(n)\le Af(Bn)+Cn$.

\begin{defn}\label{O}
A group $G$ is {\it hyperbolic relative to a collection of
subgroups} $\Hl $ if $G$ is finitely presented relative to $\Hl $
and the corresponding relative Dehn function is linear. In
particular, a group is hyperbolic (in the ordinary non--relative
sense) if and only if it is hyperbolic relative to the trivial
subgroup.
\end{defn}

The next theorem allows to regard Definition \ref{O} as a
generalization of previously known approaches. For details
concerning Bowditch's and Farb's definition of relative
hyperbolicity we refer the reader to \cite{Bow,Farb,Osin03}.

\begin{thm}[\cite{Osin03}, Theorem 1.7] \label{MainTh}
Let $G$ be a finitely generated group,  $\{ H_1, \ldots , H_m \}$
a collection of subgroups of $G$. Then the following conditions
are equivalent.

1) $G$ is finitely presented with respect to $\{ H_1, \ldots , H_m
\}$ and the corresponding relative Dehn function is linear.

2) $G$ is hyperbolic with respect to the collection $\{ H_1,
\ldots , H_m \}$ in the sense of Farb and satisfies the Bounded
Coset Penetration property (or, equivalently, $G$ is hyperbolic
relative to $\{ H_1,\ldots , H_m \}$ in the sense of Bowditch).
\end{thm}

The set of groups which have a relatively hyperbolic structure
includes fundamental groups of finite--volume non--compact
Riemannian manifolds of pinched negative curvature, geometrically
finite Kleinian groups, word hyperbolic groups, small cancellation
quotients of free products, and many other examples. This paper
continues the investigation initiated in \cite{Osin03}. It is the
second article in the sequence of three and is supposed to
establish a background for \cite{Osin04}, where we use relative
hyperbolicity to prove certain embedding theorems for countable
groups.

\begin{defn}
Let $G$ be a group hyperbolic relative to a collection of
subgroups $\{ H_\lambda , \lambda \in \Lambda \} $. A subgroup
$Q\le G$ is said to be {\it hyperbolically embedded} into $G$, if
$G$ is hyperbolic relative to $\{ H_\lambda ,\lambda \in \Lambda
\} \cup \{ Q\} $.
\end{defn}

For every element $g\in G$, we denote by $|g|_{X\cup \mathcal H}$
its {\it relative length} that is the word length with respect to
the generating set $X\cup \mathcal H$. Our main result is the
following characterization of hyperbolically embedded subgroups.

\begin{thm}\label{Q}
Suppose that $G$ is a group hyperbolic relative to a collection of
subgroups $\Hl $ and $Q$ is a subgroup of $G$. Then $Q$ is
hyperbolically embedded into $G$ if and only if the following
conditions hold:
\begin{enumerate}
\item[(Q1)] $Q$ is generated by a finite set $Y$.

\item[(Q2)] There exist $\lambda, c\ge 0$ such that for any
element $q\in Q$, we have $|q|_Y\le \lambda |q|_{X\cup \mathcal
H}+c$, where $|q|_Y$ is the word length of $q$ with respect to the
generating set $Y$ of $Q$.

\item[(Q3)] For any $g\in G$ such that $g\notin Q$, we have
$\sharp (Q\cap Q^g)< \infty $.
\end{enumerate}
\end{thm}

For ordinary hyperbolic groups the 'if' part of this theorem was
proved in \cite{Bow} (a weaker result is also obtained in
\cite{Ger}). Let us mention some corollaries of Theorem \ref{Q}.

\begin{cor}\label{Qhyp}
If $Q$ is hyperbolically embedded into $G$, then $Q$ is a
hyperbolic group.
\end{cor}

Recall that a group is called {\it elementary} if it contains a
cyclic subgroup of finite index. We also say that an element $g\in
G$ is {\it parabolic} if it is conjugate to an element of
$H_\lambda $ for some $\lambda \in \Lambda $. Otherwise $g$ is
said to be {\it hyperbolic}. In Section 3 we notify that any
hyperbolic element $g\in G$ of infinite order is contained in a
unique maximal elementary subgroup of $G$, which is denoted by
$E(g)$. Using results about cyclic subgroups of relatively
hyperbolic groups proved in \cite{Osin03}, we obtain the following
corollary of Theorem \ref{Q}.

\begin{cor}\label{Eg}
For any hyperbolic element $g\in G$ of infinite order, $E(g)$ is
hyperbolically embedded into $G$.
\end{cor}

This result can be applied to the study of boundedly generated
relatively hyperbolic groups. A group $G$ is said to be {\it
boundedly generated}, if there are elements $x_1, \ldots , x_n$ of
$G$ such that for any $g\in G$ there exist integers $\alpha _1,
\ldots , \alpha _n$ satisfying the equality $g=x_1^{\alpha
_1}\ldots x_n^{\alpha _n}$.

Bounded generation is closely related to the Congruence Subgroup
Property of arithmetic groups \cite{Rap}. It is also interesting
in connection with subgroup growth \cite{LubSeg}, unitary
representations and Kazhdan Property (T) of discrete groups
\cite{BHV,Sh}. Many lattices in semi--simple Lie groups of
$\mathbb R$--rank at least $2$ are known to be boundedly
generated. For instance, Carter and Keller \cite{CK} established
bounded generation for $SL_n(\mathcal O)$, where $\mathcal O$ is
the ring of integers of a number field and $n\ge 3$ (see also
\cite{AMe} for an elementary proof in case $\mathcal O =\mathbb Z$
and \cite{ER,K,T} for other results).

In contrast if $\Gamma $ is a uniform lattice and $G$ has $\mathbb
R$--rank $1$, then $\Gamma $ is not boundedly generated. Indeed
any such a group $\Gamma $ is non--elementary hyperbolic. The
absence of bounded generation property for a non--elementary
hyperbolic group immediately follows from the existence of
infinite periodic quotients \cite{Ols} (the direct proof can be
found in \cite{Min}). On the other hand, the problem of whether
non--uniform lattices in simple Lie groups of $\mathbb R$--rank
$1$ are boundedly generated was open until now. It is well known
that any such a lattice is hyperbolic relative to maximal
parabolic subgroups \cite{Farb,Bow95}. Thus the following theorem
answers this question negatively.

\begin{defn}
Let $G$ be a group hyperbolic relative to a collection of
subgroups $\Hl $. We say that $G$ is {\it properly hyperbolic
relative to $\Hl $}, if each of the subgroups $H_\lambda $ is
proper.
\end{defn}

\begin{thm}\label{BG}
Let $G$ be a non--elementary group properly hyperbolic relative to
a collection of subgroups $\Hl $. Then $G$ is not boundedly
generated.
\end{thm}

\section{Preliminaries}

\paragraph{Hyperbolic spaces}
Recall that a metric space $M$ is called {\it hyperbolic } (or,
more precisely,  {\it $\delta $--hyperbolic}) if for any geodesic
triangle, each side of the triangle belongs to the union of the
closed $\delta $--neighborhoods of the other two sides.

For a path $p$ in a metric space $M$, we denote by $p_-$ and $p_+$
the initial and the terminal points of $p$ respectively. All paths
under consideration are assumed to be rectifiable (i.e., of finite
length). The length of a path $p$ is denoted by $l(p)$. A path $p$
in a metric space $M$ is called $(\lambda , c)$--{\it
quasi--geodesic} for some $\lambda \ge 1$, $c\ge 0$, if for any
subpath $q$ of $p$, we have
$$l(p)\le \lambda dist_M (q_-,q_+)+c.$$
In our paper we will often use the following property of
quasi--geodesic paths in hyperbolic spaces (see \cite{Gro} or
\cite{GH}).

\begin{lem}\label{qg0}
For any $\delta \ge 0$, $\lambda \ge 1$, $c\ge 0$, there exists
$H=H(\delta , \lambda , c)$ such that for any $\delta
$--hyperbolic space, any two $(\lambda , c)$--quasi--geodesic
paths $p$, $q$ such that $p_-=q_-$, $p_+=q_+$ are contained in the
closed $H$--neighborhoods of each other.
\end{lem}

Two paths $p,q$ in a metric space $M$ are called {\it
$k$--connected}, if $$\max \{ dist_M (p_-,q_-),\; dist_M
(p_+,q_+)\} \le k.$$ The next lemma can easily be derived from the
definition of a hyperbolic space by drawing the diagonal.

\begin{lem}\label{geod}
Suppose that $p$, $q$ are $k$--connected geodesic paths in a
$\delta $--hyperbolic space and $u$ is a point on $p$ such that
$$\min \{ dist_M(u,p_-), \; dist _M(u,p_+)\} \ge 2\delta +k.$$ Then
there exists a point $v$ on $q$ such that $dist_M(u,v)\le 2\delta
$.
\end{lem}

From Lemma \ref{qg0} and Lemma \ref{geod}, we immediately obtain

\begin{cor}\label{qg}
Suppose that $p, q$ are $k$--connected $(\lambda ,
c)$--quasi--geodesic paths in a $\delta $--hyperbolic space and
$u$ is a point on $p$ such that $$\min \{ dist_M(u,p_-), \; dist
_M(u,p_+)\} \ge H+2\delta +k.$$ Then there exists a point $v$ on
$q$ such that $dist_M(u,v)\le 2(H+\delta )$.
\end{cor}

The next lemma is a simplification of Lemma 10 from \cite{Ols92}.

\begin{lem}\label{N1N2N3}
Suppose that the set of all sides of a geodesic $n$--gon
$P=p_1p_2\ldots p_n$ in a $\delta $--hyperbolic space is
partitioned into two subsets $R$ and $S$. Let $\rho $
(respectively $\sigma $) denote the sum of lengths of sides from
$R$ (respectively $S$). Assume, in addition, that $\sigma > \max
\{ \xi n,\, 10^3\rho \} $ for some $\xi \ge 3\delta\cdot 10^4$.
Then there exist two distinct sides $p_i, p_j\in S$ that contain
$13\delta $--connected segments of length greater than
$10^{-3}\xi$.
\end{lem}

\paragraph{Relatively hyperbolic groups.}
We begin with a necessary conditions for relative Dehn functions
to be well--defined.

\begin{lem}[\cite{Osin03}, Theorem 1.6]\label{DFWD}
Let $G$ be a group, $\Hl $ a collection of subgroups of $G$.
Suppose that $G$ is finitely presented with respect to $\Hl $ and
the Denh function of $G$ with respect to $\Hl $ is finite for all
values of the argument. Then the following conditions hold.

1) For any $g\in G$, the intersection $H_\lambda ^{g} \cap H_\mu $
is finite whenever $\lambda \ne \mu $.

2) The intersection $H_\lambda ^g \cap H_\lambda $ is finite for
any $g\not\in H_\lambda $.
\end{lem}

Let $G$ be a group generated by a (not necessarily finite) set
$\mathcal A$. Recall that the {\it Cayley graph} $\Ga$ of a group
$G$ with respect to the set of generators $\mathcal A$ is an
oriented labelled 1--complex with the vertex set $V(\Ga )=G$ and
the edge set $E(\Ga )=G\times \mathcal A$. An edge $e=(g,a)$ goes
from the vertex $g$ to the vertex $ga$ and has the label $\phi
(e)=a$. As usual, we denote the origin and the terminus of the
edge $e$, i.e., the vertices $g$ and $ga$, by $e_-$ and $e_+$
respectively. One can regard $\Ga $ as a metric space assuming the
length of any edge to be equal to $1$ and taking the corresponding
path metric. Also, it is easy to see that a word $W$ in $\mathcal
A$ represents $1$ in $G$ if and only if some (or, equivalently,
any) path in $\Ga $ labelled $W$ is a cycle.

Given a combinatorial path $p=e_1e_2\ldots e_k$ in the Cayley
graph $\Ga $, where $e_1, e_2, \ldots , e_k\in E(\Ga )$, we denote
by $\phi (p)$ its label. By definition, $\phi (p)=\phi
(e_1)\phi(e_2)\ldots \phi (e_k).$ We also denote by $p_-=(e_1)_-$
and $p_+=(e_k)_+$ the origin and the terminus of $p$ respectively.
A path $p$ is called {\it irreducible} if it contains no subpaths
of type $ee^{-1}$ for $e\in E(\Ga )$. The length $l(p)$ of $p$ is,
by definition, the number of edges of $p$.

In the next three lemmas we suppose that $G$ is a group hyperbolic
relative to a collection of subgroups $\Hl $. The next lemma
follows, for example, from \cite[Corollary 2.54]{Osin03}.

\begin{lem}\label{CG}
The Cayley graph $\G $ of $G$ with respect to the generating set
$X\cup \mathcal H$ is a hyperbolic metric space.
\end{lem}

\begin{lem}[\cite{Osin03}, Theorem 1.16] \label{IO}
For any hyperbolic element $g\in G$ of infinite order, there exist
positive constants $\lambda ,c$ such that
$$|g^n|_{X\cup\mathcal H}>\lambda |n|-c$$ for any $n\in \mathbb N$.
\end{lem}

\begin{lem}[\cite{Osin03}, Corollary 1.17] \label{BS}
If $g\in G$ is hyperbolic and $f^{-1}g^mf=g^n$ for some $f\in G$,
then $m=\pm n$.
\end{lem}

\paragraph{$H_\lambda $--components.}
We are going to recall an auxiliary terminology introduced in
\cite{Osin03}, which plays an important role in our paper. As
usual, by a cyclic word $W$ we mean the set of all cyclic shifts
of $W$. A word $V$ is a subword of a cyclic word $W$ if $V$ is a
subword of some cyclic shift of $W$.

\begin{defn}
Given a word $W$ (cyclic or not) in the alphabet $X\cup \mathcal
H$, we say that a subword $V$ of $W$ is an {\it $H_\lambda
$--syllable} if $V$ consists of letters from $H_\lambda\setminus\{
1\} $ and is not contained in a bigger subword entirely consisting
of letters from $H_\lambda \setminus\{ 1\} $. Let $q$ be a path
(respectively cyclic path) in $\G $. A subpath $p$ of $q$ is
called an {\it $H_\lambda $--component}, if the label of $p$ is an
$H_\lambda $--syllable of the the word $\phi (q)$ (respectively
cyclic word $\phi (q)$).
\end{defn}

\begin{defn}
Two $H_\lambda $--components $p_1, p_2$ of a path $q$ (cyclic or
not) in $\G $ are called {\it connected} if there exists a path
$c$ in $\G $ that connects some vertex of $p_1$ to some vertex of
$p_2$ and ${\phi (c)}$ is a word consisting of letters from $
H_\lambda\setminus\{ 1\} $. In algebraic terms this means that all
vertices of $p_1$ and $p_2$ belong to the same coset $gH_\lambda $
for some $g\in G$. Note that we can always assume that $c$ has
length at most $1$, as every element of $H_\lambda $ is included
in the set of generators.  An $H_\lambda $--component $p$ of a
path $q$ (cyclic or not) is called {\it isolated } if no
(distinct) $H_\lambda $--component of $q$ is connected to $p$.
\end{defn}

The next lemma is a simplification of Lemma 2.27 from
\cite{Osin03}. The subsets $\Omega _\lambda $ mentioned below are
exactly the sets of all elements of $H_\lambda $ represented by
$H_\lambda $--components of defining words $R\in \mathcal R$ in a
suitably chosen finite relative presentation $\langle X,\;
H_\lambda, \lambda\in \Lambda \; |\; R=1,\; R\in \mathcal R
\rangle $ of $G$.

\begin{lem}\label{Omega}
Suppose that $G$ is a group hyperbolic relative to a collection of
subgroups $\Hl $. Then there exists a constant $K>0$ and subsets
$\Omega _\lambda \subseteq H_\lambda $ such that the following
conditions hold.

1) The union $\Omega =\bigcup\limits_{\lambda\in \Lambda }\Omega
_\lambda $ is finite.

2) Let $q$ be a cycle in $\G $, $p_1, \ldots , p_k$ a set of
isolated $H_\lambda $--components of $q$ for some $\lambda\in
\Lambda $, $g_1, \ldots , g_k$ the elements of $G$ represented by
the labels of $p_1, \ldots , p_k$ respectively. Then for any $i=1,
\ldots , k$, $g_i$ belongs to the subgroup $\langle \Omega
_\lambda \rangle \le G$ and the lengths of $g_i$ with respect to
$\Omega _\lambda $ satisfy the inequality $$ \sum\limits_{i=1}^k
|g_i|_{\Omega _\lambda }\le Kl(q).$$
\end{lem}

\section{Hyperbolically embedded subgroups}

Throughout this section we fix a group $G$ hyperbolic relative to
a collection of subgroups $\{ H_\lambda ,\; \lambda \in \Lambda
\}$, a finite relative generating set $X=X^{-1}$ of $G$ with
respect to $\{ H_\lambda ,\; \lambda \in \Lambda \}$, and the set
$\Omega $ provided by Lemma \ref{Omega}. By $\dxh $ we denote the
natural metric on the Cayley graph $\G $.

\begin{lem}\label{Q1}
For any $\lambda  \ge 1$, $c\ge 0$, there exists $\alpha_1=\alpha
_1(\lambda , c)>0$ such that for any $k\ge 0$ there exists $\alpha
_2=\alpha _2(k, \lambda , c)>0$ satisfying the following
condition. Let $p$, $q$ be two $k$--connected $(\lambda ,
c)$--quasi--geodesics in $\G$, such that the labels of $p$ and $q$
are words in the alphabet $X$. Let $u$ be a vertex on $p$ such
that $$\min\{ \dxh (u,p_-),\; \dxh (u, p_+)\} \ge \alpha _2.$$
Then there exists a vertex $v$ on $q$ such that the element
$u^{-1}v$ belongs to the subgroup $\langle X\cup \Omega \rangle $
and the length of $u^{-1}v$ with respect to $X\cup \Omega $
satisfies the inequality $|u^{-1}v|_{X\cup \Omega }\le \alpha_1$.
\end{lem}

\begin{proof}
By Lemma \ref{CG}, the Cayley graph $\G $ is a hyperbolic metric
space with respect to the metric $\dxh $. Let $\delta $ denote the
hyperbolicity constant of $\G $. We set
$$\alpha _2= 5H+6\delta +k ,$$ where $H=H(\delta , \lambda , c)$
is the constant provided by Lemma \ref{qg0}.

\begin{figure}
\begin{picture}(100,36.87)(0,60)
\qbezier[2000](16.97,93.87)(81.14,80.43)(102.18,93.87)
\qbezier[2000](16.79,66.29)(80.96,79.73)(102,66.29)
\put(16.62,94.04){\circle*{.79}} \put(42.78,89.09){\circle*{.79}}
\put(42.78,70.89){\circle*{.79}} \put(16.44,66.11){\circle*{.79}}
\put(46.06,79.36){\circle*{.79}} \put(66.94,76.49){\circle*{.79}}
\put(67.06,83.36){\circle*{.79}} \put(102.18,93.87){\circle*{.79}}
\put(102,66.29){\circle*{.79}} \put(68.06,87.15){\circle*{.79}}
\put(67.88,73.01){\circle*{.79}}
\qbezier(67.88,87.15)(65.41,77.96)(67.88,73.01)
\qbezier(42.78,89.09)(49.5,81.58)(42.78,70.89)
\put(85.56,69.83){\makebox(0,0)[cc]{$p$}}
\put(84.85,90.16){\makebox(0,0)[cc]{$q$}}
\put(15.91,64.17){\makebox(0,0)[cc]{$p_-$}}
\put(15.85,96){\makebox(0,0)[cc]{$q_-$}}
\put(68.06,70.89){\makebox(0,0)[cc]{$u$}}
\put(67.88,89.8){\makebox(0,0)[cc]{$v$}}
\put(43,68.4){\makebox(0,0)[cc]{$u_0$}}
\put(41.72,91.75){\makebox(0,0)[cc]{$v_0$}}
\put(105,95){\makebox(0,0)[cc]{$q_+$}}
\put(104.6,65){\makebox(0,0)[cc]{$p_+$}}

\put(45.62,84){\vector(-1,4){.07}}

\put(66.83,81){\vector(0,1){.07}}

\put(54,79.8){\vector(1,0){.07}}

\put(48,83.62){\makebox(0,0)[cc]{$s$}}
\put(65,81){\makebox(0,0)[cc]{$t$}}
\put(55.5,77){\makebox(0,0)[cc]{$e$}}
\put(43,79.25){\makebox(0,0)[cc]{$s_-$}}
\put(70,76.5){\makebox(0,0)[cc]{$t_-$}}

\put(85.62,88.25){\vector(4,1){.07}}

\put(85.87,71.87){\vector(4,-1){.07}}

\put(45,86){\circle*{.79}}

\thicklines \qbezier(44.87,86)(46.25,82.37)(46.12,79.25)
\qbezier(67,83.5)(66.37,78.94)(66.75,76.62)
\qbezier(46.12,79.25)(58.31,81.12)(66.75,76.5)
\end{picture}
\caption{}\label{Q-fig1}
\end{figure}

Without loss of generality we may assume that $4(H+\delta )$ is
integer. Since $\dxh (p_-,u)\ge \alpha _2$, there exists a vertex
$u_0$ on the segment $[p_-, u]$ such that
\begin{equation}\label{Q11}
\dxh (u_0,u)=4(H+\delta )
\end{equation}
Note that $$\max\{\dxh (p_-,u_0),\; \dxh (p_+, u_0)\} \ge
H+2\delta +k.$$ Hence by Corollary \ref{qg}, there exist points
$v_0$ and $v$ on $q$ such that
\begin{equation}\label{Q12}
\dxh (u,v)\le 2(H+\delta ),
\end{equation}
and
\begin{equation}\label{Q121}
\dxh (u_0,v_0)\le 2(H+\delta ).
\end{equation}

Clearly we may assume that $v$ and $v_0$ are vertices of $p$. Let
us consider the combinatorial loop
$$r=[u,u_0][u_0,v_0][v_0,v][u,v]^{-1},$$ where $[u_0,v_0]$
and $[u,v]$ are arbitrary geodesics in $\G $ and $[u,u_0]$,
$[v_0,v]$ are segments of $p^{-1}$ and $q$ (or $q^{-1}$)
respectively. Obviously inequalities (\ref{Q11}), (\ref{Q12}) and
(\ref{Q121}) imply
$$
\begin{array}{rl}
l([v_0,v])\le & \lambda \dxh (v_0,v)+c\le \\ & \\ & \lambda (\dxh
(v_0, u_0)+\dxh (u_0,u)+\dxh (u,v))+c\le \\ & \\ & 8\lambda
(H+\delta )+c.
\end{array}
$$
Therefore,
\begin{equation}\label{Q13}
\begin{array}{rl}
l(r) \le & \lambda \dxh (u,u_0)+c +\dxh (u_0,v_0)+l([v_0,v])+\dxh
(v,u)\le \\ & \\ & (12\lambda +4)(H+\delta )+2c.
\end{array}
\end{equation}

We are going to show that for any $\lambda \in \Lambda$, any
$H_\lambda $--component of $[u,v]$ is isolated in $r$. Indeed
assume that a certain $H_\lambda $--component $t$ of $[u,v]$ is
not isolated in $r$. Since the labels $\phi ([u_0,u])$ and $\phi
([v_0,v])$ are words in $X$, they contain no $H_\lambda
$--components at all. Thus the only possibility is that there
exists an $H_\lambda $--component $s$ of $[u_0, v_0]$ that is
connected to $t$. This means that there exists a path $e$ of
length at most $1$ in $\G $ connecting $s_-$ to $t_-$ (see Fig.
\ref{Q-fig1}). Using (\ref{Q12}) and (\ref{Q121}) we obtain
$$
\begin{array}{rl}
\dxh (u_0,u)\le & \dxh (u_0, s_-)+\dxh (s_-, t_-)+\dxh (t_-, u)\le
\\ & \\ &(\dxh (u_0, v_0)-1) +1+(\dxh (v, u)-1)< \\ & \\ & 4(H+\delta )
\end{array}
$$
which contradicts to (\ref{Q11}).

Thus for any $\lambda \in \Lambda$, any $H_\lambda $--component
$t$ of $[u,v]$ is isolated in $r$. Therefore, by Lemma \ref{Omega}
we have $t\in \langle X\cup \Omega \rangle$ and $|t|_{X\cup \Omega
}\le Kl(r)$, where $K$ depends on $G$ only. Combining these with
(\ref{Q13}) and (\ref{Q12}), we obtain $u^{-1}v \in \langle X\cup
\Omega \rangle $ and
\begin{equation}\label{uv}
|u^{-1}v|_{X\cup \Omega } \le Kl(r)l([u,v])\le 2K(H+\delta )(
(12\lambda +4)(H+\delta )+2c) .
\end{equation}
It remains to set $\alpha _1$ to be equal to the right side of
(\ref{uv}).
\end{proof}

In the next two lemmas $Q$ is a subgroup of $G$ satisfying
conditions (Q1)--(Q3). Since $\sharp Y<\infty $, without loss of
generality we may assume that $Y\subseteq X$.

\begin{lem}\label{Q2}
For every $\alpha > 0$ there exists $A=A(\alpha )>0$ such that the
following holds. If $a,b\in Q$ and $f,g\in G$ are arbitrary
elements such that $\max \{ |a|_Y, |b|_Y\} \ge A$, $\max \{
|f|_{X\cup\mathcal H}, |g|_{X\cup\mathcal H}\} \le \alpha $, and
$a=fbg$, then $f,g\in Q$.
\end{lem}

\begin{proof}
Let $\alpha _1=\alpha _1(\lambda ,c)$ and $\alpha _2=\alpha _2
(\alpha , \lambda , c)$ be constants provided by Lemma \ref{Q1},
where $\lambda $ and $c$ are given by (Q2). We also denote by $N$
the number of different elements of the subgroup $\langle X\cup
\Omega \rangle $ of length at most $\alpha _1$ with respect to the
generating set $X\cup \Omega $. Since $\sharp (X\cup \Omega
)<\infty $, $N<\infty $. By (Q3), there is an integer $M>0$ such
that for any element $t\in \langle X\cup \Omega \rangle $ of
length $|t|_{X\cup \Omega}\le \alpha _1$, any element of $Q\cap
Q^t$ has length strictly less than $M$ with respect to $Y$
whenever $t\notin Q$. Set
$$A=2(\lambda \alpha _2+ c)+MN.$$

Increasing constants if necessary, we can assume that $\lambda $,
$c$, $\alpha _2$, and $M$ are integer. By the condition of the
lemma there is a quadrangle $rqsp^{-1}$ in $\G $, where $r$ and
$s$ are geodesics in $\G $ whose labels represent $f$ and $g$, and
$p$, $q$ are paths labelled by the shortest words in $Y\subseteq
X$ representing $a$ and $b$ respectively. Note that paths $p$ and
$q$ are $(\lambda ,c)$--quasi--geodesic according to (Q2). We take
a vertex $a_0$ on $p$ such that the length of the segment
$[p_-,a_0]$ of $p$ equals $\lambda \alpha _2 +c$. By the choice of
$A$ the length of the segment $[a_0,p_+]$ of $p$ is at least
$\lambda \alpha _2+c+MN$. Let $a_0, \ldots , a_N$ denote the
subsequent vertices of $[a_0,p_+]$ such that the length of the
segment $[a_{k-1},a_k]$ of $p$ is equal to $M$ for each $k=1,
\ldots , N$. Note that for any $i=1, \ldots , N$,
$$\min \{ \dxh (p_-, a_i),\,\dxh (a_i,p_+)\} \ge
 \frac{\min\{ l([p_-,a_i]),\, l([a_i,p_+])\} -c}{\lambda } \ge
\alpha _2. $$

\begin{figure}
\begin{picture}(95.28,35)(-7,45)
\qbezier(13.08,77.96)(62.84,69.03)(91.04,76.72)
\qbezier(12.9,55.33)(62.67,64.26)(90.86,56.57)
\put(60.63,73.15){\vector(1,0){.07}}

\put(60.46,60.1){\vector(1,0){.07}}

\put(39.6,73.89){\circle*{1}} \put(13.26,77.96){\circle*{1}}
\put(13.08,55.33){\circle*{1}} \put(91.22,56.57){\circle*{1}}
\put(91.04,76.36){\circle*{1}}
\put(82.21,58.16){\circle*{1}} \put(41.01,59.04){\circle*{1}}
\put(24.04,57.1){\circle*{1}} \put(66.11,59.75){\circle*{1}}
\put(65.05,73.01){\circle*{1}}
\put(70.71,56.1){\circle*{.5}} \put(29.52,55.57){\circle*{.5}}
\put(47.55,56.45){\circle*{.5}}
\put(73.72,56.1){\circle*{.5}} \put(32.53,55.57){\circle*{.5}}
\put(50.56,56.45){\circle*{.5}}
\put(76.72,56.1){\circle*{.5}} \put(35.53,55.57){\circle*{.5}}
\put(53.56,56.45){\circle*{.5}}
\put(56.04,62.22){\makebox(0,0)[cc]{$p$}}
\put(55.86,71.42){\makebox(0,0)[cc]{$q$}}


\multiput(39.42,74.07)(.0331454,-.3130395){48}{\line(0,-1){.3130395}}

\multiput(64.88,73.18)(.033145,-.414317){32}{\line(0,-1){.414317}}


\put(11.67,53.21){\makebox(0,0)[cc]{$p_-$}}
\put(23.86,54.8){\makebox(0,0)[cc]{$a_0$}}
\put(65.94,56.6){\makebox(0,0)[cc]{$a_j$}}
\put(40.66,55.51){\makebox(0,0)[cc]{$a_i$}}
\put(82.91,55.51){\makebox(0,0)[cc]{$a_N$}}
\put(41,76.72){\makebox(0,0)[cc]{$b_i$}}
\put(64.35,76.19){\makebox(0,0)[cc]{$b_j$}}
\put(93.5,54){\makebox(0,0)[cc]{$p_+$}}
\put(11.31,80.61){\makebox(0,0)[cc]{$q_-$}}
\put(92.98,79.2){\makebox(0,0)[cc]{$q_+$}}

\put(40.13,66.64){\vector(-1,4){.07}}

\put(65.41,67){\vector(0,1){.07}}

\put(37,67){\makebox(0,0)[cc]{$t_i$}}
\put(68.2,67.53){\makebox(0,0)[cc]{$t_j$}}
\qbezier(13.08,77.96)(18.22,67.44)(12.9,55.15)
\qbezier(91.22,76.54)(86.13,66.03)(91.04,56.57)
\put(15.7,67){\vector(0,1){.07}}
\put(88.66,65.6){\vector(0,-1){.07}} \put(13,67){$r$}
\put(90,67){$s$}

\end{picture}

  \caption{}\label{Q-fig2}
\end{figure}

Since $l(r)=|f|_{X\cup\mathcal H}$ and $l(s)=|g|_{X\cup\mathcal
H}$, the paths $p$ and $q$ are $\alpha $--connected. By Lemma
\ref{Q1}, there are vertices $b_0, \ldots , b_N$ on $q$ and paths
$t_0, \ldots t_N$ such that $(t_k)_-=a_k$, $(t_k)_+=b_k$, and
$\phi (t_k)$ represents an element of the subgroup $\langle X\cup
\Omega \rangle $ of length at most $\alpha _1$ with respect to
$X\cup \Omega $ (see Fig. \ref{Q-fig2}). By our choice of the
constant $N$, there are two paths $t_i$, $t_j$ such that $\phi
(t_i)$ and $\phi (t_j)$ represent the same element $t$ in $G$.
Reading the label of the cycle
$t_i[b_i,b_j]t_j^{-1}[a_i,a_j]^{-1}$, where $[b_i,b_j]$ is the
segment of $q$, gives us the equality $tq_1t^{-1}=q_2$ for some
$q_1,q_2\in Q$. By the choice of vertices $a_0,\ldots , a_N$, we
have $|q_2|_Y=l([a_i, a_j])\ge M$. Therefore, $t\in Q$ according
to the choice of $M$. Since $\phi ([p_-,a_i]t_i[q_-,b_i]^{-1})$
represents the same element of $G$ as $\phi (r)$ and $\phi ([p_-,
a_i])$, $\phi ([q_-, b_i])$ represent elements of $Q$, we obtain
$f\in Q$. Similarly $g\in Q$.
\end{proof}

Recall that $G$ is generated by $X$ relative to $\Hl $, $Q$ is
generated by $Y$, and $Y\subseteq X$. Let $Z=X\setminus Y$.  We
consider the groups
$$ F=\left( \ast _{\lambda\in \Lambda } H_\lambda
\right) \ast F(Y)\ast F(Z) $$ and
$$F_Q=\left( \ast _{\lambda\in
\Lambda } H_\lambda  \right) \ast Q\ast F(Z)$$ (here $F(Y)$ and
$F(Z)$ stand for free groups freely generated by $Y$ and $Z$
respectively) together with the following commutative diagram of
homomorphisms defined in the obvious way.

\unitlength 0.7mm
\begin{picture}(28.28,30)(-62,2)
\put(25.99,24.22){\vector(1,0){.07}}
\put(7.42,24.22){\line(1,0){18.562}}
\put(17.32,9.37){\vector(-3,-4){.07}}
\multiput(27.75,21.92)(-.033644597,-.040487566){310}{\line(0,-1){.040487566}}
\put(15.56,9.9){\vector(3,-4){.07}}
\multiput(5.83,21.74)(.033642624,-.040982832){289}{\line(0,-1){.040982832}}
\put(3.71,24.4){\makebox(0,0)[cc]{$F$}}
\put(32,24.22){\makebox(0,0)[cc]{$F_Q$}}
\put(16.6,6){\makebox(0,0)[cc]{$G$}}
\put(16.79,26.87){\makebox(0,0)[cc]{$\e $}}
\put(6.72,15.2){\makebox(0,0)[cc]{$\beta $}}
\put(26.5,14.85){\makebox(0,0)[cc]{$\gamma $}}
\end{picture}

Since $G$ is hyperbolic with respect to $\{ H_\lambda ,\lambda \in
\Lambda \},$ there is a finite subset $\mathcal R$ such that
$Ker\, \beta=\langle \mathcal R\rangle ^F$. Then obviously $Ker\,
\gamma =\langle \e (R),\; R\in\mathcal R\rangle ^{F_Q}$. In
particular, $G$ is finitely presented with respect to the
collection $\Hl \cup \{ Q\} $. To simplify our notation we denote
by $\mathcal Q$ alphabet $X\cup \mathcal H\sqcup (Q\setminus \{
1\} )$. For any word $T$ in the alphabet $X\cup \mathcal H$ such
that $T\in Ker\, \beta$, we denote by $Area^{rel}(T)$ the minimal
possible number $k$ of factors in the decomposition
\begin{equation}\label{prod1}
T=_{F}\prod\limits_{i=1}^k f_i^{-1}R_if_i,
\end{equation}
where $f_i\in F$ and $R_i\in \mathcal R$, $i=1, \ldots , k$.
Similarly, given a word $W$ in $\mathcal Q$ such that $W\in Ker\,
\gamma $, we denote by $Area^{rel}_Q(W)$ the minimal possible
number $k$ of factors in the decomposition
\begin{equation}\label{prod1}
W=_{F_Q}\prod\limits_{i=1}^k g_i^{-1}\e (R_i)g_i,
\end{equation}
where $g_i\in F_Q$ and $R_i\in \mathcal R$, $i=1, \ldots , k$. The
next lemma is quite obvious and we leave the proof to the reader.

\begin{lem} \label{relar}
Let $U,V,W$ be words in $\mathcal Q$, $T$ a word in $X\cup\mathcal
H$. Suppose that $U,V\in Ker\, \gamma $ and $T\in Ker\, \beta $.
Then:
\begin{enumerate}
\item[(a)] $Area ^{rel}_Q(UV)\le Area^{rel}_Q(U)+Area^{rel}_Q(V).$

\item[(b)] $Area ^{rel}_Q(W^{-1}VW)=Area^{rel}_Q(V).$

\item[(c)] $Area^{rel}_Q(\e (T))\le Area^{rel}(T)$.
\end{enumerate}
\end{lem}

In the proof of Theorem \ref{Q}, we will use the following
auxiliary notion.

\begin{defn}
A word $W$ in $\mathcal Q$ is called {\it primitive} if $W$ can
not be decomposed as $$W\equiv W_1q_1W_2q_2W_3,$$ where
$q_1,q_2\in Q\setminus \{ 1\} $ and the subword $W_2$ represents
an element of the subgroup $Q$ in $G$.
\end{defn}

\begin{lem}\label{prim}
Suppose that there exists a constant $\kappa >0$ such that for any
primitive word $W$ in $\mathcal Q$, $W\in Ker\, \gamma $, we have
$Area ^{rel}_{F_Q} (W)\le \kappa \| W\| $. Then the relative Dehn
function of $G$ with respect to $\Hl \cup \mathcal \{ Q\} $ is
linear.
\end{lem}

\begin{proof}
Let $U$ be an arbitrary word in $\mathcal Q$ such that $U\in Ker\,
\gamma $. To prove the lemma it suffices to show that
\begin{equation} \label{U}
Area^{rel}_{F_Q} (U)\le \kappa \| U\| .
\end{equation}
Let $q(U)$ denote the number of letters from $Q\setminus \{ 1\} $
that appear in $U$. We proceed by induction on $q(U)$. If
$q(U)=0$, the inequality (\ref{U}) obviously holds since $U$ is
primitive in this case.

Now let $q$ be non--primitive and $q(U)>0$. Then $U\equiv
U_1q_1U_2q_2U_3$, where $q_1,q_2\in Q\setminus \{ 1\} $ and
$\gamma (U_2)\in Q$. Let $r_1,r_2$ be letters from the alphabet
$Q\setminus \{ 1\} $ that represent the same elements as $\gamma
(U_2)$ and $\gamma (q_1U_2q_2)$ in $G$ respectively. In
particular, we have $\gamma (r_2)=_{G}\gamma (q_1r_1q_2)$.
However, the restriction of $\gamma $ to $Q$ is injective. Hence
$r_2=_{F_Q}q_1r_1q_2$.  Therefore,
$$
\begin{array}{rl}
U=_{F_Q} &
(U_1q_1r_1q_2U_3)(U_3^{-1}q_2^{-1}r_1^{-1}U_2q_2U_3)=_{F_Q}\\ &\\
& (U_1r_2U_3)(U_3^{-1}q_2^{-1}r_1^{-1}U_2q_2U_3).
\end{array}
$$
Note that $q(U_1r_2U_3)<q(U)$ and $q(r_1^{-1}U_2)<q(U)$. Moreover,
\begin{equation}\label{qU}
q(U_1r_2U_3)+q(r_1^{-1}U_2) = q(U_1)+1+q(U_3)+1+q(U_2)=q(U).
\end{equation}
Applying Lemma \ref{relar} together with inequality (\ref{qU}) and
the inductive hypothesis, we obtain
$$
\begin{array}{cl}
Area^{rel}_{Q} (U)\le & Area^{rel}_{Q}(U_1r_2U_3)
+Area^{rel}_{Q}(r_1^{-1}U_2)\le \\ & \\
& \kappa \| U_1r_2U_3\| +\kappa \| r_1^{-1}U_2\| \le \kappa \|
U\|.
\end{array}
$$
\end{proof}

Now we are ready to prove the main theorem.

\begin{proof}[Proof of Theorem \ref{Q}]
First suppose that $Q$ satisfies (Q1)--(Q3). We keep the notation
introduces above. In particular, we assume that $Y$ is a subset of
$X$.

Since $G$ is hyperbolic relative to $\Hl $, there exists $L>0$
such that for any word $V$ in $X\cup\mathcal H$ such that $V\in
Ker\, \beta $, we have $Area ^{rel}(V)\le Ln$. Let us take a word
$W$ in $\mathcal Q$, $W\in Ker\, \gamma $, of length $\| W\|=n$.
We want to bound $Area^{rel}_Q (W)$ from above by a linear
function of $n$. Taking into account Lemma \ref{prim}, we may
assume that $W$ is primitive. In case $W$ contains no letters from
$Q\setminus \{ 1\} $, we immediately obtain
\begin{equation}\label{triv}
Area^{rel}_Q(W)\le Ln.
\end{equation}
In what follows we assume that at least one letter from
$Q\setminus \{ 1\} $ appears in $W$.

Let $$W\equiv W_0t_1W_1\ldots t_lW_l,$$ where $t_i\in Q\setminus
\{ 1\} $ and subwords $W_0, \ldots , W_l$ contain no letters from
$Q\setminus \{ 1\} $. For each $t_i$, $i=1, \ldots l$, we fix a
shortest word $V_i$ in the alphabet $Y$ such that $V_i=_Q t_i$ and
consider the word
$$V\equiv W_0V_1W_1\ldots V_lW_l$$ in the alphabet $X\cup
\mathcal H$ regarded as an element of the group $F$. Clearly
$\epsilon (V)=W$. We set
$$\theta=\sum\limits_{i=1}^l \| V_i\|,$$
$$\omega =\sum\limits_{i=1}^l \| W_i\|.$$
Note that
\begin{equation}\label{on}
\omega <n.
\end{equation}

Also let $\lambda $, $c$ be given by (Q2), $\delta $ the
hyperbolicity constant of $\G $, $H=H(\delta , \lambda , c)$ the
constant from Lemma \ref{qg0}, $A=A(13\delta )$ the constant
provided by Lemma \ref{Q2}, and
$$\xi =\max \{ 3\delta \cdot 10^4,\; (A+2H)\cdot 10^3\} .$$ There
are three possibilities to consider.

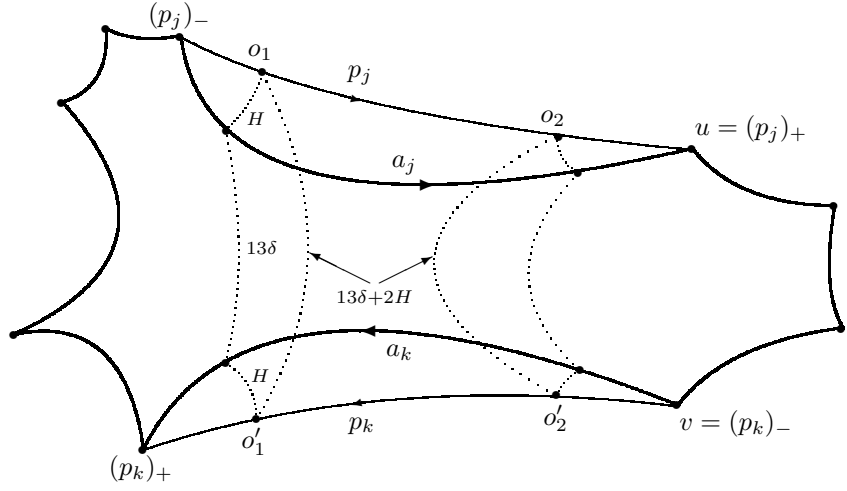
\begin{figure}
\unitlength 1mm 
\linethickness{0.4pt}
\ifx\plotpoint\undefined\newsavebox{\plotpoint}\fi 
\begin{picture}(126.34,62.5)(10,5)

\qbezier(38.01,64.7)(57.36,54.18)(106.07,49.67)
\qbezier(33.06,9.72)(64.35,20.24)(104.12,15.56)
\qbezier[50](48.97,59.75)(61.43,35.36)(48.08,13.79)
\qbezier[50](88.03,16.97)(55.6,34.74)(87.86,51.44)

\put(28.11,65.76){\circle*{1}} \put(22.23,55.89){\circle*{1}}
\put(38.01,64.7){\circle*{1}} \put(48.97,60.1){\circle*{1}}
\put(88.39,51.44){\circle*{1}} \put(106.07,49.85){\circle*{1}}
\put(124.94,42.23){\circle*{1}} \put(125.94,25.85){\circle*{1}}
\put(104.07,15.73){\circle*{1}} \put(91.19,20.35){\circle*{1}}
\put(88.07,16.98){\circle*{1}} \put(90.94,46.6){\circle*{1}}
\put(48.19,13.85){\circle*{1}} \put(44.07,21.35){\circle*{1}}
\put(44.19,52.23){\circle*{1}} \put(33.07,9.73){\circle*{1}}
\put(15.94,24.98){\circle*{1}} \put(68,45){\line(1,0){.5}}
\put(60.88,15.75){\vector(-4,-1){.07}}
\put(62.13,56.25){\vector(4,-1){.07}}

\put(38,67.5){\makebox(0,0)[cc]{$(p_j)_-$}}
\put(114,52.5){\makebox(0,0)[cc]{$u=(p_j)_+$}}
\put(112,13.25){\makebox(0,0)[cc]{$v=(p_k)_-$}}
\put(32.88,7.25){\makebox(0,0)[cc]{$(p_k)_+$}}
\put(48.75,62.63){\makebox(0,0)[cc]{$o_1$}}
\put(87.63,53.88){\makebox(0,0)[cc]{$o_2$}}
\put(88,14){\makebox(0,0)[cc]{$o_2^\prime $}}
\put(47.88,10.88){\makebox(0,0)[cc]{$o_1^\prime $}}
\put(67.88,47.63){\makebox(0,0)[cc]{$a_j$}}
\put(62,59.63){\makebox(0,0)[cc]{$p_j$}}
\put(67.25,22.75){\makebox(0,0)[cc]{$a_k$}}
\put(62.25,13){\makebox(0,0)[cc]{$p_k$}}
\put(49,36.63){\makebox(0,0)[cc]{$_{13\delta }$}}
\put(48.2,53.75){\makebox(0,0)[cc]{$_H$}}
\put(48.8,19.63){\makebox(0,0)[cc]{$_H$}}
\put(64,30){\makebox(0,0)[cc]{$_{13\delta +2H}$}}

\put(64.6,32){\vector(2,1){7}} \put(63.3,32){\vector(-2,1){8}}

\qbezier[35](44.02,52.5)(47.55,34.83)(44.02,21.39)
\qbezier[15](44.02,52.33)(47.2,54.98)(48.97,59.75)

\qbezier[13](44.02,21.39)(48.35,18.47)(48.08,13.79)
\qbezier[10](88.03,51.62)(88.3,47.73)(91.04,46.67)
\qbezier[30](91.04,46.67)(77.6,32.97)(91.04,20.33)
\qbezier[7](91.04,20.33)(89.01,19.18)(88.03,16.97)

\thicklines

\qbezier(38.01,64.88)(41.81,35.27)(106.07,49.85)
\qbezier(106.07,49.85)(111.81,42.34)(124.98,42.25)
\qbezier(124.98,42.25)(123.04,31.29)(126.04,25.99)
\qbezier(126.04,25.99)(110.49,24.04)(104.12,15.73)
\qbezier(104.12,15.73)(47.38,38.1)(33.06,9.9)
\qbezier(33.06,9.9)(30.51,28.22)(15.78,25.1)
\qbezier(15.78,25.1)(40.04,36.15)(22.27,56.04)
\qbezier(22.27,56.04)(28.55,57.89)(28.11,65.76)
\qbezier(28.11,65.76)(32.35,63.73)(38.01,64.88)
\put(72,45.05){\vector(1,0){.07}}
\put(62,25.57){\vector(-1,0){.07}}

\end{picture}
\caption{$(13\delta +2H)$--connected segments of the paths $p_j$
and $p_k$.}\label{41}
\end{figure}

a) First assume that $\theta \le (\lambda \xi +c)n$. Taking into
account inequality (\ref{on}), we obtain
$$ \| V\| =\theta+\omega < (\lambda \xi +c+1)n.$$ By Lemma
\ref{relar}, we have
\begin{equation} \label{A1}
Area^{rel}_Q (W)\le Area^{rel} (V)< L(\lambda \xi +c+1)n.
\end{equation}

b) Further suppose that $\theta \le 10^3\lambda \omega +cn$.
Similarly to the previous case, we obtain
$$\| V\|\le 10^3\lambda \omega +cn+\omega < (10^3\lambda +c+1)n $$ and
\begin{equation} \label{A2}
Area^{rel}_Q (W)\le Area^{rel} (V)< L(10^3\lambda +c+1)n.
\end{equation}

c) Finally, suppose that
\begin{equation}\label{casec}
\theta >\max\{ (\lambda \xi +c)n,\, 10^3\lambda \omega +cn\} .
\end{equation}
We consider a cycle $q$ in the Cayley graph $\G $ labelled $V$.
Let $q=q_0p_1q_1\ldots p_lq_l$, where $q_i$, $p_i$ are subpaths
labelled $W_i$ and $V_i$ respectively. Consider the $2l+1$--gon
$b_0a_1b_1\ldots a_lb_l$, where $a_i$ (respectively $b_i$) is a
geodesic path with the same endpoints as $p_i$ (respectively
$q_i$). Let $S=\{ a_i, i=1, \ldots , l\} $, $R=\{ b_i, i=0, \ldots
, l\} $. By $\rho $ (respectively $\sigma $) we denote the sum of
lengths of sides from $R$ (respectively $S$). Obviously $\rho \le
\omega $. Further, according to (Q2) $p_i$ is $(\lambda ,
c)$--quasi--geodesic in $\G $ for any $i=1, \ldots , l$. In
particular, we have
$$l(a_i)\ge \frac1\lambda (\| V_i\| -c).$$ Hence,
$$
\sigma \ge \sum\limits_{i=1}^l \frac1\lambda (\| V_i\| -c)=
\frac1\lambda (\theta -cl)\ge \frac1\lambda (\theta -cn).
$$
Together with (\ref{casec}), this yields
$$
\sigma > \max\{ \xi n,\, 10^3\omega \} \ge \max\{ \xi n,\,
10^3\rho \} .
$$

By Lemma \ref{N1N2N3}, there are two sides, say $a_j$ and $a_k$,
$k>j$, having $13\delta $--bounded segments of length at least
$\xi \cdot 10^{-3}$. Therefore, by Lemma \ref{qg0}, $p_j$ and
$p_k$ have $(13\delta +2H)$--bounded segments $s=[o_1,o_2]$,
$s^\prime=[o_1^\prime , o_2^\prime ]$ of length at least $\xi
\cdot 10^{-3}-2H\ge A$ (see Fig. \ref{41}). Set $a=o_1^{-1}o_2$,
$f=o_1^{-1}o_1^\prime $, $b=(o_1^\prime )^{-1} o_2^\prime $,
$g=(o_2^\prime )^{-1}o_2$. Obviously the elements $a,b,f,g $
satisfy the requirements of Lemma \ref{Q2}. Thus $f,g\in Q$. Let
$(p_j)_+=u$, $(p_k)_-=v$. Then
$$u^{-1}v=(u^{-1}o_1)f^{-1}((o^\prime _1)^{-1}v)\in Q,$$ as labels
of the segments $[u,o_1]$ and $[o_1, v]$ of $p$ and $q$
respectively are words in $Y$ and hence represent elements of $Q$.

Clearly the words $$W_0t_1\ldots W_{j-1}t_j$$ and $$W_0t_1\ldots
W_{k-2}t_{k-1}W_{k-1}$$ represent $u$ and $v$ respectively in $G$.
Therefore, the subword $$W_jt_{j+1}\ldots W_{k-2}t_{k-1}W_{k-1}$$
of $W$ represents the element $u^{-1}v\in Q$. However this
contradicts to the assumption that $W$ is primitive. Thus case c)
is impossible for primitive $W$.

Taking together (\ref{triv}), (\ref{A1}), and (\ref{A2}), we
obtain  $$Area^{rel}_Q (W)< \kappa n,$$ where $$\kappa =\max\{
L(10^3\lambda +c+1), L(\lambda \xi +c+1)\}.$$ This completes the
proof of the first part of the theorem.

Now suppose that the subgroup $Q$ is hyperbolically embedded into
$G$. By $\Gamma (G, \mathcal Q)$ we denote the Cayley graph of $G$
with respect to $\mathcal Q$. Also let $\Omega _Q$ denote the
subset of $Q$ given by Lemma \ref{Omega} applied to the collection
of subgroups $\Hl\cup \{ Q\} $.

For every nontrivial element $q\in Q$, we fix a shortest word $W$
over $X\cup \mathcal H$ that represents $q$ in $G$. Let us
consider a cycle $d=pr$ in $\Gamma (G, \mathcal Q)$, where $p$ is
an edge labelled $q$ and $\phi (r)\equiv W^{-1}$. Clearly $p$ is
an isolated $Q $--component of $q$ as $r$ contains no edges
labelled by elements of $Q$. Applying Lemma \ref{Omega}, we obtain
$q\in \langle \Omega _Q \rangle $. Therefore $\Omega _Q $
generates $Q$. Moreover, we have
$$|q|_{\Omega _Q }\le K_Ql(d)\le K_Q (\| W\| +1) \le
K_Q(|q|_{X\cup \mathcal H}+1),$$ where $K_Q$ is some constant
independent of $q$. Thus the conditions (Q1) and (Q2) hold. The
fulfilment of (Q3) follows from Lemma \ref{DFWD}.
\end{proof}

Recall that a group is hyperbolic if it is finitely generated and
its Cayley graph with respect to some finite generating set is a
hyperbolic metric space.

\begin{proof}[Proof of Corollary \ref{Qhyp}.]
Let $Y$ be a finite generating set of $Q$. As above we assume that
$Y$ is a subset of $X$. The inclusion $Y\subseteq X$ defines the
embedding of the Cayley graph $$\iota : \Gamma (Q,Y)\to \G $$
whose restriction on $Q=V(\Gamma (Q,Y))$ is the identity map. If
$\Delta $ is a geodesic triangle in $\Gamma (Q,Y)$, then $\iota
(\Delta )$ is a triangle in $\G $ whose sides are $(\lambda ,
c)$--quasi--geodesic according to (Q2). Let $v$ be a vertex on a
side of $\Delta $. If $\delta $ is the hyperbolicity constant of
$\G $, then there is a vertex $w$ on the union of the other two
sides of $\Delta $ such that
$$|v^{-1}w|_{X\cup\mathcal H}=\dxh (\iota(v),\iota (w))\le
2H+\delta ,$$ where $H=H(\delta , \lambda , c)$ is the constant
from Lemma \ref{qg0}. Therefore, by (Q2) the distance between $v$
and $w$ in $\Gamma (Q,Y)$ is $$|v^{-1}w|_Y\le \lambda |v^{-1}w|
+c\le \lambda (2H+\delta )+c.$$ This shows that $\Gamma (Q,Y)$ is
$(\lambda (2H+\delta ) +c+1)$--hyperbolic.
\end{proof}

\section{Elementary subgroups and bounded generation}

All assumptions and notation listed at the beginning of the
previous section remain valid here. We begin with auxiliary
results about elementary subgroups of relatively hyperbolic
groups.

\begin{lem}\label{E1}
For any hyperbolic element of infinite order $g\in G$, there
exists a constant $C=C(g)$ such that if $f^{-1}g^nf=g^n$ for some
$f\in G$ and some $n\in \mathbb N$, then there are $m\in \mathbb
Z$ and $h\in \langle X\cup\Omega \rangle $ such that $f=hg^m$ and
$|h|_{X\cup \Omega }\le C$.
\end{lem}

\begin{proof}
Without loss of generality we may assume that $g\in X$. Let
$\lambda , c$ be constants from Lemma \ref{IO} and $\alpha
_1=\alpha _1(\lambda , c)$, $\alpha _2=\alpha _2(k,\lambda , c)$
be constants provided by Lemma \ref{Q1}, where $k=|f|_{X\cup
\mathcal H}$.

Since $fg^n=g^nf$, we have $fg^{nt}=g^{nt}f$ for any $t\in \mathbb
Z$. Consider the $k$--connected paths $p$, $q$ in $\G $ labelled
$g^{nt}$ such that $q_-=1$, $q_+=g^{nt}$, $p_-=f$,
$p_+=fg^{nt}=g^{nt}f$. If $t$ is big enough, there is a vertex $u$
on $p$ such that $u=fg^{nj}$ for some $j\in \mathbb Z$ and $$\min
\{ \dxh (p_-, u),\; \dxh(p_+,u)\} \ge \alpha _2.$$ By Lemma
\ref{Q1}, there exists a vertex $v$ on $q$ such that
$|v^{-1}u|_{X\cup\Omega }\le \alpha _1$. Note that $v=g^i$ for
some $i\in \mathbb Z $. We consider the vertex $w=g^{nj}$ on $q$.
The equality $fg^{nj}=g^{nj}f$ implies $ f=w^{-1} u =(w^{-1}v)
(v^{-1}u)=g^{i-nj}h$ for $h=v^{-1}u$.
\end{proof}

\begin{defn}
For any hyperbolic element of infinite order $g\in G$, we set
$$E(g)=\{ f\in G\; :\; f^{-1}g^nf=g^{\pm n}\; {\rm for \; some\; } n\in
\mathbb N\} .$$
\end{defn}

\begin{thm}
Every hyperbolic element $g\in G$ is contained in a unique maximal
elementary subgroup, namely in $E(g)$.
\end{thm}

\begin{proof}
Let $$E_+(g)=\{ f\in G\; :\; f^{-1}g^nf=g^n\; {\rm for \; some\; }
n\in \mathbb N\} .$$ By Lemma \ref{E1}, $\langle g\rangle $ has
finite index in $E_+(g)$. Clearly $E_+(g)$ has index $2$ in
$E(g)$. Hence $|E(g):\langle g\rangle |< \infty $. In particular,
$E(g)$ is elementary.

It remains to show that if $E$ is another elementary subgroup
containing $g$ then $E\le E(g)$. Let $s$ be an element of $E$ such
that $\langle s\rangle $ has finite index in $E$. Passing to a
subgroup $\langle s^i\rangle $ for some $i$ if necessary, we may
assume that $\langle s\rangle $ is normal in $E$. Obviously
$s^l=g^k$ for some $k,l\in \mathbb Z\setminus \{ 0\} $. In
particular, $s$ is hyperbolic. Indeed if $s\in H_\lambda ^a$ for
some $\lambda \in \Lambda $, $a\in G$, then the cyclic subgroup
$\langle s^l\rangle $ is contained in the intersection $H_\lambda
^{ag}\cap H_\lambda ^a=\left( H_\lambda ^{aga^{-1}}\cap H_\lambda
\right) ^a $ as $s^l$ commutes with $g$. Therefore, $H_\lambda
^{aga^{-1}}\cap H_\lambda $ is infinite. By Lemma \ref{DFWD},
$aga^{-1}\in H_\lambda $ that contradicts to hyperbolicity of $g$.
Thus $s$ is hyperbolic and for any element $t\in E$, we have
$t^{-1}st=s^{\pm 1} $ by Corollary \ref{BS}. Hence
$$t^{-1}g^kt=t^{-1}s^lt=s^{\pm l}=g^{\pm k}.$$ By the definition
of $E(g)$, we have $t\in E(g)$.
\end{proof}

\begin{proof}[Proof of Corollary \ref{Eg}.]
Let us show that conditions (Q1)--(Q3) hold for $Q=E(g)$. The
fulfillment of (Q1) is obvious. Let $Y$ be a finite generating set
of $E(g)$ containing $g$. Since $\langle g\rangle $ has finite
index in $E(g)$, there is a constant $D>0$ such that any element
$t\in E(g)$ can be represented as $t=hg^m$ for some $m\in \mathbb
Z$, where $\max \{ |h|_{X\cup \mathcal H} ,\; |h|_Y\} \le D$.
Using Lemma \ref{IO} we obtain
$$ |t|_Y\le |g^m|_Y+D\le |m|+D\le \frac{1}{\lambda } |g^m|_{X\cup \mathcal
H}+c+D\le \frac1\lambda \left( |t|_{X\cup \mathcal H}+D\right)
+c+D $$ for some positive $\lambda $, $c$. Thus (Q2) holds.

Finally if the intersection $E(g)^f\cap E(g)$ is infinite for some
$f\in G$, then it contains $g^n$ for some $n\in \mathbb Z\setminus
\{ 0\} $. Hence $f^{-1}g^nf=g^k$ for some $k\in \mathbb Z\setminus
\{ 0\} $. By Corollary \ref{BS} $k=\pm n$ and thus $f\in E(g)$.
\end{proof}

\begin{lem}\label{ah}
For any $\lambda \in \Lambda $ there is a finite subset $\mathcal
F_\lambda \subseteq H_\lambda $ such that if $h\in H_\lambda
\setminus \mathcal F_\lambda $, $a\in G\setminus H_\lambda $, and
$|a|_{X\cup\mathcal H}=1$, then $ah$ is a hyperbolic element of
infinite order.
\end{lem}

\begin{proof} Set $$\mathcal F_\lambda =\{ f\in \langle \Omega
_\lambda \rangle ,\; |f|_{\Omega \lambda }\le 5K\} ,$$ where $K$
and $\Omega _\lambda $ are given by Lemma \ref{Omega}. Since $a\in
G\setminus H_\lambda $ has relative length $1$, we can think of
$a$ is a letter from $X\cup \mathcal H$. For every $m\in \mathbb
N$, we consider a path $p_m=q_1r_1q_2r_2\ldots q_mr_m$, where
$q_i$ (respectively $r_i$) is labelled $a$ (respectively $h$),
$i=1, \ldots , m$.

\begin{figure}
\unitlength=1mm \linethickness{0.4pt}
\begin{picture}(140,23)(-25,4)
\qbezier(65.23,7.95)(36.15,35.89)(6.01,7.95)

\put(5.83,7.95){\circle*{1}} \put(15.94,15.68){\circle*{1}}
\put(55.92,15.23){\circle*{1}} \put(50.28,18.65){\circle*{1}}
\put(22.03,18.95){\circle*{1}} \put(65.29,7.8){\circle*{1}}
\put(37.01,21.9){\vector(1,0){.07}}

\put(13,9){\makebox(0,0)[cc]{$p_m$}}
\put(16.5,19.62){\makebox(0,0)[cc]{$r_i$}}
\put(54.5,19.3){\makebox(0,0)[cc]{$r_j$}}
\put(35,11.4){\makebox(0,0)[cc]{$e$}}
\put(35.64,25){\makebox(0,0)[cc]{$s$}}

\thicklines \linethickness{1pt}
\qbezier(15.76,15.61)(18.28,16.95)(22,18.88)
\qbezier(50.09,18.58)(53.74,16.65)(55.89,15.31)
\qbezier(22,18.88)(33.52,9.51)(50.09,18.58)
\put(54.07,16.6){\vector(2,-1){.07}}
\put(20,17.7){\vector(2,1){.07}}
\put(36.7,14.1){\vector(1,0){.07}}

\end{picture}
  \caption{}\label{ahm}
\end{figure}

Note that $r_1, \ldots , r_m$ are $H_\lambda $--components of
$p_m$. First of all we are going to show that they are isolated.
Indeed suppose that $r_i$ is connected to $r_j$ for some $j>i$ and
$j-i$ is minimal possible. Let $s$ denote the segment $[(r_i)_+,
(r_j)_-]$ of $p_m$, and let $e$ be a path of length at most $1$ in
$\G $ labelled by an element of $H_\lambda $ such that
$e_-=(r_i)_+$, $e_+=(r_j)_-$ (see Fig. \ref{ahm}). If $j=i+1$,
then $\phi (s)=a $. Since $\phi (s)$ and $\phi (e)$ represent the
same element in $G$, we arrive at a contradiction with the
assumption $a\notin H_\lambda $. Therefore, $j=i+1+k$ for some
$k\ge 1.$ Note that the components $r_{i+1}, \ldots , r_{i+k}$ are
isolated in the cycle $se^{-1}$. (Otherwise we can pass to another
pair of connected $H_\lambda $--components with smaller value of
$j-i$.) By Lemma \ref{Omega} we have $h\in \langle \Omega _\lambda
\rangle $ and
$$
k |h|_{\Omega _\lambda }\le Kl(se^{-1})=K (2k+2).
$$
Hence $|h|_{\Omega _\lambda }\le K (2+2/k)\le 4K$ which
contradicts to the choice of $h$. Thus all components $r_1, \ldots
, r_m $ are isolated in $p_m$. In particular, this means that the
element $ah$ has infinite order. Indeed if $(ah)^n=1$ for some
$n\in \mathbb N$, then the components $r_1$ and $r_{n+1}$ of
$p_{n+1}$ coincide (and thus they are connected).

Let us show that $ah$ is hyperbolic. Indeed suppose that
$ah=b^{-1}fb$ for some $b\in G$, $f\in H_\nu $, $\nu \in \Lambda
$. We take $m=4|b|_{X\cup\mathcal H}+6$. Let $B$ be a shortest
word in $X\cup \mathcal H$ representing $b$ in $G$. Consider a
cycle $d=p_mq_m$ in $\G $ such that $\phi (q_m)\equiv
(B^{-1}tB)^{-1}$, where $t$ is the letter from $H_\nu\setminus \{
1\} $ representing the same element as $f^m$ in $G$. Note that
$r_2, \ldots , r_{m-1}$ are components of $d$ and any $H_\lambda
$--component of the subpath $q_m$ of $d$ is connected to at most
one $H_\lambda $--component from the set $\{ r_2, \ldots ,
r_m{m-1}\} $. (If an $H_\lambda $--component of $q_m$ is connected
to $r_i$ and $r_j$, then $r_i$ and $r_j$ are connected.) As the
total number of $H_\lambda $--components in $q_m$ does not exceed
$l(q_m)\le 2|b|_{X\cup\mathcal H}+1$, at least
$(m-2|b|_{X\cup\mathcal H}-3)=m/2$ components from the set $\{
r_2, \ldots , r_{m-1} \} $ are isolated. Lemma \ref{Omega} yields
$$
\frac{1}{2}m |h|_{\Omega _\lambda } \le Kl(d)\le K\left( 2m
+2|b|_{X\cup\mathcal H}+1\right) <\frac{5}{2}mK.
$$
Dividing by $m/2$, we obtain $|h|_{\Omega _\lambda }\le 5K$. This
contradicts to $h\notin \mathcal F_\lambda $. The lemma is proved.
\end{proof}

\begin{cor} \label{prop}
Let $G$ be an infinite group properly hyperbolic relative to a
collection of subgroups $\Hl $. Then $G$ contains a hyperbolic
element of infinite order.
\end{cor}

\begin{proof}
Removing trivial subgroups from the set $\Hl $ if necessary, we
may assume that $H_\lambda \ne \{ 1\} $ for any $\lambda \in
\Lambda $. First suppose that $\sharp\, \Lambda =\infty $. Since
$G$ is defined by a finite presentation (\ref{G}) with respect to
$\Hl $, there is a finite subset $\Lambda _0\subseteq \Lambda $
such that no relators from $\mathcal R$ involve letters from
$H_\lambda $ for $\lambda \in \Lambda \setminus \Lambda _0$.
Therefore $G=G_0\ast \left( \ast_{\lambda \in \Lambda
\setminus\Lambda _0} H_\lambda \right) $, where $G_0$ is the
subgroup of $G$ generated by all $H_\lambda $, $\lambda \in
\Lambda _0$. If $\lambda _1, \lambda _2\in\Lambda \setminus\Lambda
_0$, $\lambda _1\ne \lambda _2$, then for any two nontrivial
elements $g_1\in H_{\lambda _1}$, $g_2\in H_{\lambda _2} $, the
product $g_1g_2$ is hyperbolic and has infinite order.

Now suppose that $\sharp \, \Lambda <\infty $. If all subgroups
$H_\lambda $ are finite, then $\G $ is locally finite. As $\G $ is
hyperbolic by Lemma \ref{CG}, $G$ is a hyperbolic group (in the
ordinary non--relative sense). It is well--known that any infinite
hyperbolic group contains an element of infinite order $g$.
Obviously $g$ is hyperbolic in this case, as all parabolic
elements of $G$ have finite orders.

Finally if there is an infinite subgroup in the collection $\{
H_\lambda , \lambda \in \Lambda \} $, the desired hyperbolic
element exists by Lemma \ref{ah}. The element $a\in G\setminus
H_\lambda $ of relative length $1$ mentioned in Lemma \ref{ah}
exists as $G$ is generated by elements of relative length $1$ and
$H_\lambda \ne G$.
\end{proof}

Corollary \ref{prop} together with Lemma \ref{IO} immediately
imply the following.

\begin{cor} \label{InfD}
If $G$ is infinite and properly hyperbolic relative to a
collection of subgroups $\Hl $, then $G$ has infinite diameter
with respect to the metric $\dxh $.
\end{cor}

Finally let us prove Theorem \ref{BG}.

\begin{proof}[Proof of Theorem \ref{BG}]
We assume that $G$ is non--elementary and is properly hyperbolic
relative to $\Hl $.  Suppose that there are elements $x_1, \ldots
, x_n$ of $G$ such that for any $g\in G$ there exist integers
$\alpha _1, \ldots , \alpha _n$ satisfying the equality
$g=x_1^{\alpha _1}\ldots x_n^{\alpha _n}$. According to Corollary
\ref{Eg}, if $x_i$ is a hyperbolic element of infinite order for
some $i$, then $G$ is properly hyperbolic relative to the
collection $\Hl \cup\{ E(x_i)\} $. In case $x_i$ has finite order,
$G$ is properly hyperbolic relative to the collection $\Hl \cup\{
\langle x_i\rangle \} $ by Theorem \ref{Q}, as any finite subgroup
satisfies (Q1)--(Q3). Thus joining new subgroups to the collection
$\Hl $ if necessary, we may assume that elements $x_1, \ldots ,
x_n$ are hyperbolic, i.e., for any $i=1, \ldots , n$,
$x_i=a_i^{-1}h_ia_i$, where $h_i\in H_{\lambda _i}$ for some
$\lambda _i\in \Lambda $. Then for any integers $\alpha _1, \ldots
, \alpha _n$, we have
$$x_1^{\alpha _1}\ldots x_n^{\alpha _n}=a_1h_1^{\alpha _1}a_1^{-1}
\ldots a_nh_n^{\alpha _n}a_n^{-1} .$$ Therefore, for every element
$g=x_1^{\alpha _1}\ldots x_n^{\alpha _n}$, we have
$$|g|_{X\cup \mathcal H}\le \sum\limits_{i=1}^n (2|a_i|_{X\cup
\mathcal H} +1).$$ This means that $G $ has finite diameter with
respect to $\dxh $ contradictory to Corollary \ref{InfD}.
\end{proof}

\vspace{1cm}

\noindent Department of Mathematics\\ Vanderbilt University\\
Nashville,
TN 37240\\

\noindent {\it E-mail address:} denis.ossine@vanderbilt.edu

\end{document}